\documentclass[a4paper, BCOR=0.0mm, DIV=calc]{scrartcl}
\usepackage[T1]{fontenc}
\usepackage[utf8]{inputenc}
\usepackage[ngerman]{babel}
\usepackage[osf,sc]{mathpazo}
\usepackage{textcomp}
\usepackage{microtype}
\usepackage{amsmath, amssymb, amsfonts}
\KOMAoptions{twoside=false, twocolumn=false, headinclude=false, footinclude=false, mpinclude=false, pagesize=auto}
\recalctypearea
\newcommand{\ebinom}[2]{\left(\frac{#1}{#2} \right)}
\begin{document}
\boldmath
\title{Beweis einer einzigartigen Transformation von Reihen\footnote{
Originaltitel: "`Specimen transformationis singularis serierum"', erstmals publiziert in "`\textit{Nova Acta Academiae Scientarum Imperialis Petropolitinae} 12, 1801, pp. 58-70"', Nachdruck in "`\textit{Opera Omnia}: Series 1, Volume 16, pp. 41 - 55 "', Eneström-Nummer E710, übersetzt von: Alexander Aycock, Textsatz: Artur Diener,  im Rahmen des Projektes "`Eulerkreis Mainz"' }}
\unboldmath

\author{Leonhard Euler}
\date{}
\maketitle
\paragraph{§1}
Ich habe diese Reihe betrachtet
\[
	s = 1 + \frac{ab}{1\cdot c}x + \prod \frac{(a+1)(b+1)}{2\cdot(c+1)}x^2 + \prod \frac{(a+2)(b+2)}{3\cdot (c+2)}x^3 + \mathrm{etc.}
\]
wo auf gewohnte Weise $\Pi$ den Koeffizienten des vorhergehenden Terms bezeichnet. Diese Reihe ist so beschaffen, dass ihre Summe im Allgemeinen auf keine Weise ausgedrückt werden zu können scheint, obwohl sie trotzdem in allen Fällen, in denen entweder $a$ oder $b$ eine ganze negative Zahl ist, abbricht und ihre Summe auf endliche Weise ausgedrückt wird.

\paragraph{§2}
Wenn wir nun also
\[
	s = z(1-x)^{c-a-b}
\]
setzen und weiter
\[
	c - a = \alpha ~~ \mathrm{und} ~~ c - b = \beta
\]
wird der Buchstabe $z$ die Summe dieser der vorhergehenden ähnlichen 
\[
	z = 1 + \frac{\alpha\beta}{1\cdot c}x + \prod \frac{(\alpha + 1)(\beta + 1)}{2(\alpha + 1)}x^2 + \prod \frac{(\alpha + 2)(\beta + 2)}{3(c+2)}x^3 + \mathrm{etc.}
\]
ausdrücken, die nun in allen Fällen abbricht, in denen entweder $\alpha$ oder $\beta$ eine negative Zahl ist, und daher sooft entweder $a-c$ oder $b-c$ eine positive Zahl war.

\paragraph{§3}
Diese Transformation ist von umso größerer Bedeutung anzusehen, weil sie außer durch lange Um-- und Irrwege und sogar durch Differentialgleichungen zweiter Ordnung gefunden werden zu können scheint. Daher wird es der Mühen wert sein, die ganze Analysis, auf der diese Transformation beruht, klar erörtert zu haben.

\paragraph{§4}
Weil
\[
	s = 1 + \frac{ab}{1\cdot c}x + \frac{ab}{1\cdot c}\cdot \frac{(a+1)(b+1)}{2(c+1)}\cdot xx + \mathrm{etc.}
\] 
ist und sogar in jedem folgenden Term der Zähler wie der Nenner zwei neue Faktoren erhält, wollen wir durch Differentiation zuerst aus jedem Term die zwei letzten Faktoren entfernen, weil durch diese Operatorenen
\[
	\frac{\partial s}{\partial x} = \frac{ab}{1\cdot c} + \frac{ab}{1\cdot c}\frac{(a+1)(b+1)}{c+1}x + \mathrm{etc.}
\]
dargestellt werden wird, welche Gleichung mit $x^c$ multipliziert und erneut differenziert 
\[
	\partial \cdot x^c \partial s = abx^{c-1} + \frac{ab}{1\cdot c}(a+1)(b+1)x^c + \mathrm{etc.}
\]
liefert, wo wir, die wir immer nach Kürze suchen, das Element $\partial x$ weggelassen haben, was man sich selbst ein wenig merken kann.

\paragraph{§5}
Gleich wollen wir auf ähnliche Weise durch Differentiation den einzelnen Zählern zwei neue Faktoren hinzufügen, und zwar auf diese Weise:
\begin{enumerate}
	\item Unsere Reihe wird mit $x^a$ multipliziert und differenziert
	\[
		\partial\cdot x^a s = ax^{a-1} + \frac{ab}{1\cdot c}(a+1)x^a + \mathrm{etc}
	\] 
	ergeben, welche
	\item mit $x^{b+1-a}$ multipliziert und wiederum differenziert
	\[
		\partial\cdot x^{b+1-a}\partial x^a s = abx^{b-1} + \frac{ab}{1\cdot c}(a+1)(b+1)x^b + \mathrm{etc.}
	\]
\end{enumerate}
liefert, welche Form aus der vorhergehenden entsteht, wenn sie mit $x^{b-c}$ multipliziert wird.

\paragraph{§6}
Daher erhalten wir also diese Gleichung
\[
	\partial\cdot x^{b-a+1}\partial\cdot x^a s = x^{b-c}\partial\cdot x^c \partial s
\]
welche entwickelt auf diese Form
\[
	x^{b+1}\partial\partial s + (a+b+1)x^b\partial s + abx^{b-1}s = x^b\partial\partial s + cx^{b-1}\partial s
\]
zurückgeführt wird. Diese Gleichung nimmt durch $x^{b-1}$ geteilt, und nachdem alle Terme auf die rechte Seite gebracht wurden, diese Form an
\[
	0 = x(1-x)\partial\partial s + \left[ c - (a+b+1)x \right]\partial s - abs
\]
so dass von der Auflösung dieser Differentialgleichung zweiter Ordnung die Summation der vorgelegten Reihe abhängt. Aber diese Gleichung scheint in der Tat so beschaffen zu sein, dass sie im Allgemeinen keine Integration zulässt.

\paragraph{§7}
Obwohl aber diese Differentialgleichung uns wenig an Hilfe zu bringen scheint, ist sie trotzdem für eine ausgezeichnete Transformation empfänglich, mit welcher man unsere ganze Aufgabe erledigt. Wir benutzen nämlich diese allgemeine Substitution
\[
	s = (1-x)^n z
\]
woraus
\[
	\log{s} = n\log{(1-x)} + \log{z}
\]
wird, und durch Differenzieren wird
\[
	\frac{\partial s}{s} = \frac{\partial z}{z} - \frac{n\partial x}{1-x}
\]
sein, welche Gleichung erneut differenziert
\[
	\frac{\partial\partial s}{s} - \frac{\partial s^2}{ss} = \frac{\partial\partial z}{z} - \frac{\partial z^2}{zz} - \frac{n\partial x^2}{(1-x)^2}
\]
liefert. Dieser füge man nun diese Gleichung
\[
	\frac{\partial s^2}{ss} = \frac{\partial z^2}{zz} - \frac{2n\partial x\partial z}{z(1-x)} + \frac{nn\partial x^2}{(1-x)^2}
\]
hinzu und es wird
\[
	\frac{\partial\partial s}{s} = \frac{\partial\partial z}{z} - \frac{2n\partial x\partial z}{z(1-x)} + \frac{n(n-1)\partial x^2}{(1-x)^2}
\]
hervorgehen.

\paragraph{§8}
Wenn also gleich die vorgelegte Gleichung durch $s$ geteilt so dargestellt wird
\[
	0 = x(1-x)\frac{\partial\partial s}{s} + \left[ c - (a+b+1)x \right] \frac{\partial s}{s} - abs
\]
werden wir nach Ausführung der Substitution zu einer Differentialgleichung zweiter Ordnung zwischen $z$ und $x$ gelangen, welche
\begin{align*}
	&x(1-x)\frac{\partial\partial z}{z} - \frac{2nx\partial x\partial z}{z} + \left[ c-(a+b+1)x\right] \frac{\partial z}{z} \\ 
	&+ \frac{n(n-1)x\partial x^2}{1-x} - \frac{n[c-(a+b+1)x]\partial x}{1-x} - ab = 0
\end{align*}
sein wird.

\paragraph{§9}
Es ist aber klar, dass hier die Zahl $n$ so angenommen werden kann, dass die letzten Glieder, die den Nenner $1-x$ enthalten, durch sie geteilt werden können; das passiert im Fall $n = -a-b+c$, nach Einführung welchen Wertes, sodass $s = (1-x)^{c-a-b}z$ ist, und so die Gleichung zwischen $z$ und $x$ diese Form erhalten wird:
\[
	x(1-x)\partial\partial z + \left[ c+(a+b-2c-1)x\right]\partial z - (c-a)(c-b)z = 0
\]

\paragraph{§10}
Wenn wir also gleich in dieser Gleichung $c-a = \alpha$ und $c-b = \beta$ setzen, wird die Gleichung zwischen $z$ und $x$ unter dieser Form erscheinen
\[
	x(1-x)\partial\partial z + \left[ c - (\alpha + \beta + 1)x\right]\partial z - \alpha \beta z = 0
\]
welche von der ersten überhaupt nicht abweicht, außer dass wir anstelle der Buchstaben $a$ und $b$  hier $\alpha$ und $\beta$ haben. Weil daher die erste Differenzen-Differentialgleichung aus der Reihe
\[
	s = 1 + \frac{ab}{1\cdot c}x + \prod\frac{(a+1)(b+1)}{2(c+1)}xx + \prod\frac{(a+2)(b+2)}{3(c+2)}x^3 + \mathrm{etc.}
\]
entstanden ist, wird andererseits aus der zweiten Gleichung die Reihe
\[
	z = 1 + \frac{\alpha\beta}{1\cdot c}x + \prod\frac{(\alpha + 1)(\beta + 1)}{2(c+1)}xx + \prod\frac{(\alpha + 1)(\beta + 1)}{3(c+2)}x^3 + \mathrm{etc.}
\]
entstehen, während $\alpha = c-a$ und $\beta = c-b$ ; und diese zwei Reihen $s$ und $z$ hängen so voneinander ab, dass $s = (1-x)^{c-a-b}z$ ist oder $\frac{s}{z} = (1-x)^{c-a-b}$.

\paragraph{§11}
Aber aus der ersten Differenzen-Differentialgleichung kann mit einer direkten Methode dieselbe Reihe für $z$ gefunden werden. Weil nämlich aus der ersten Reihe für $x=0$ gesetzt $s=1$ wird, wir nun aber $z = (1-x)^{a+b-c}s$ gesetzt haben, wird in demselben Fall $x=0$ gleich $z=s=1$ werden. Nachdem das bemerkt wurde, wollen wir für $z$ diese Reihe ansetzen:
\[
	z = 1 + Ax + Bx^2 + Cx^3 + Dx^4 + \mathrm{etc.}
\]
woher
\[
	\partial z = A  + 2Bx + 3Cx^2 + 4Dx^3 + 5Ex^4 + \mathrm{etc.}
\]
wird und
\[
	\partial\partial z = 2B + 6Cx + 12Dx^2 + 20Ex^3 + 30Fx^4 + \mathrm{etc.}
\]
nach Einsetzen welcher Werte \\

\begin{tabular}{rlllll}
$x(1-x)\partial\partial z$ &= & $2Bx$ & $+6Cx^2$ & $+12Dx^3$ & $\mathrm{+etc.}$ \\
& & & $-2Bx^2$ & $-6Cx^3$ & $\mathrm{-etc.}$ \\
$c\partial z$ &= $Ac$ & $+Bcx$ & $+3Ccx^2$ & $+4Dcx^3$ & $\mathrm{+etc.}$ \\
$-(\alpha + \beta + 1)x\partial z$ &= & $-(\alpha + \beta + 1)Ax$ & $-2(\alpha + \beta +1)Bx^2$ & $-3(\alpha + \beta +1)Cx^3$ & $\mathrm{-etc.}$ \\
$-\alpha\beta z$ &= $-\alpha\beta$ & $-A\alpha\beta x$ & $-B\alpha\beta x^2$ & $-C\alpha\beta x^3$ & $\mathrm{-etc.}$
\end{tabular}

\hrule
\[
	x(1-x)\partial\partial z + c\partial z - (\alpha + \beta + 1)x\partial z - \alpha\beta z = 0
\]
hervorgehen wird.

\paragraph{§12}
Nachdem also alle Terme null gesetzt wurden, werden wir die folgenden Gleichungen erhalten:\\[1em]
\begin{tabular}[c]{lcl}
&I.& $Ac - \alpha\beta = 0$ \\
&II.& $2B(c+1) - (\alpha + 1)(\beta + 1)A = 0$ \\
&III.& $3C(c+2) - (\alpha + 2)(\beta + 2)B = 0$ \\
&IV.& $4D(c+3) - (\alpha + 3)(\beta + 3)C = 0$ \\
&V.& $5E(c+4) - (\alpha + 4)(\beta + 4)D = 0$ \\
&$\mathrm{etc.}$ &
\end{tabular}

\paragraph{§13}
Daher werden also dieselben Koeffizienten gefunden, die wir schon hatten, natürlich

\begin{tabular}[c]{lcl}
&$A$ &$=$ $\dfrac{\alpha\beta}{1\cdot c}$ \\[4mm]
&$B$ &$=$ $\dfrac{A(\alpha + 1)(\beta + 1)}{2(c+1)}$ \\[4mm]
&$C$ &$=$ $\dfrac{B(\alpha + 2)(\beta + 2)}{3(c+2)}$ \\[4mm]
&$D$ &$=$ $\dfrac{C(\alpha + 3)(\beta + 3)}{4(c + 3)}$ \\[4mm]
&$\mathrm{etc.}$&
\end{tabular}\\[1em]

Weil aber die Methode, mit welcher wir diese herausragende Transformation erhalten haben, besonders skuril ist und durch lange Umwege vorgeht, wäre besonders zu wünschen, dass man eine andere direktere und natürlichere Methode fände, durch die man gewiss einen nicht zu verachtenden Zuwachs in die Analysis einbrächte. Ich gestehe aber, dass ich mich bisher bei dieser Untersuchung vergebens bemüht habe.

\paragraph{§14}
Weil in diesen Reihen die Anzahl der Faktoren stetig wächst, wählen wir, damit wir hier die aus der Binomialpotenz entstehenden Charaktere angenehmer benutzen können, die Buchstaben $a$ und $b$, ebenso wie die Buchstaben $\alpha$ und $\beta$; negative Werte teilen wir zu, indem wir $a = -f$, $b=-g$, $\alpha = - \zeta$ und $\beta = - \eta$ setzen, so dass $\zeta = -c-f$ und $\eta = -c-g$ ist, und schon werden unsere beiden Reihen $s$ und $z$ voneinander abhängen, sodass
\[
	s = (1-x)^{c+f+g}z
\]
ist. Wir wollen nun zuerst gemäß dieser Werte die erste Reihe $s$ entwickeln, und sie wird
\[
	s = 1 + \frac{fg}{1\cdot c}x + \prod\frac{(f-1)(g-1)}{2(c+1)}x^2 + \prod\frac{(f-2)(g-2)}{3(c+2)}x^3 + \mathrm{etc.}
\]
sein und auf ähnliche Weise wird die zweite Reihe
\[
	z = 1 + \frac{\zeta\eta}{1\cdot c}x + \prod\frac{(\zeta - 1)(\eta - 1)}{2(c+1)}x^2 + \prod\frac{(\zeta - 2)(\eta - 2)}{3(c+2)}x^3 + \mathrm{etc.}
\]
sein.

\paragraph{§15}
Hier können wir schon angenehm die erwähnten Charaktere anwenden. Es bezeichne also $\left(\frac{m}{n}\right)$ den Koeffizienten des Termes $v^n$, welcher mit selbigem aus der Entwicklung der Binomialpotenz $(1+v)^m$ übereinstimmt, sodass wir auf diese Weise
\[
	(1+v)^m = 1 + \left(\frac{m}{1}\right)v + \left(\frac{m}{2}\right)v^2 + \left(\frac{m}{3}\right)v^3 + \mathrm{etc.}
\]
haben. Daher wird also für die erste unserer Reihen $\frac{f}{1} = \left(\frac{f}{1}\right)$ werden; darauf
\[
	\frac{f(f-1)}{1\cdot 2} = \left(\frac{f}{2}\right); ~~ \frac{f(f-1)(f-2)}{1\cdot 2\cdot 3} = \left(\frac{f}{3}\right) ~~ \mathrm{etc.}
\]
und so wird diese Reihe gleich angenehmer so angegeben:
\[
	s = 1 + \frac{g}{c}\left(\frac{f}{1}\right)x + \frac{g}{c}\cdot\frac{g-1}{c+1}\cdot\left(\frac{f}{2}\right)x^2 + \frac{g}{c}\cdot\frac{g-1}{c+1}\cdot\frac{g-2}{c+2}\cdot\left(\frac{f}{3}\right)x^3 + \mathrm{etc.}
\]
Damit wir nun die Terme, die den Buchstaben $g$ enthalten, auf ähnliche Weise zusammenziehen, wollen wir die Reihe selbst auf beiden Seiten mit dem Charakter $\left(\frac{g+c-1}{c-1}\right)$ multiplizieren; dann wird nämlich $\left(\frac{g+c-1}{c-1}\right)\frac{g}{c} = \left(\frac{g+c-1}{c}\right)$ sein; $\left(\frac{g+c-1}{c-1}\right)\frac{g}{c}\cdot\frac{g-1}{c+1} = \left(\frac{g+c-1}{c+1}\right)$, welcher Ausdruck weiter mit $\frac{g-2}{c+2}$ multipliziert diesen Charakter geben wird: $\left(\frac{g+c-1}{c+2}\right)$. Nachdem diese Dinge bemerkt wurden, erhalten wir nun diese Reihe:
\begin{align*}
s\left(\frac{g+c-1}{c-1}\right) =& \left(\frac{g+c-1}{c-1}\right) + \left(\frac{f}{1}\right)\left(\frac{g+c-1}{c}\right) x \\
&+ \left(\frac{f}{2}\right)\left(\frac{g+c-1}{c+1}\right)x^2 + \left(\frac{f}{3}\right)\left(\frac{g+c-1}{c+2}\right) x^3 + \mathrm{etc.}
\end{align*}

\paragraph{§16}
Auf ähnliche Weise wird sich auch die andere Reihe transformieren lassen; dort ist aber natürlich zu bemerken, dass diese Transformation auf zweierlei Weise aufgestellt werden kann, je nachdem, ob die Faktoren des Nenners $1,2,3,4, \mathrm{etc}$ entweder mit dem Buchstaben $\zeta$ oder mit $\eta$ verbunden werden. Zuerst werden wir also aus der vorhergehenden Reihe, wenn wir $\zeta$ anstelle von $f$ und $\eta$ anstelle von $g$ schreiben, diese Reihe erhalten:
\begin{align*}
z\ebinom{\eta + c -1}{c -1} =& \ebinom{\eta + c -1}{c -1} + \ebinom{\zeta}{1}\ebinom{\eta +c -1}{c}x + \ebinom{\zeta}{2}\ebinom{\eta +c -1}{c+1}x^2 \\
&+ \ebinom{\zeta}{3}\ebinom{\eta +c -1}{c+2}x^3 + \mathrm{etc.}
\end{align*}
Wenn wir aber anstelle von $f$ und $g$ in umgekehrter Reihenfolge $\eta$ und $\zeta$ schreiben, geht
\begin{align*}
z\ebinom{\zeta +c -1}{c-1} =& \ebinom{\zeta +c -1}{c-1} + \ebinom{\eta}{1}\ebinom{\zeta +c -1}{c}x + \ebinom{\eta}{2}\ebinom{\zeta +c -1}{c+1}x^2 \\
&+ \ebinom{\eta}{3}\ebinom{\zeta +c -1}{c+2}x^3 + \mathrm{etc.}
\end{align*}
hervor. Für jede von beiden aber bleibt die Relation dieselbe, natürlich
\[
	s = (1-x)^{c+f+g}z
\]

\paragraph{§17}
Damit es klarer erscheint, wie sehr sich diese zwei Reihen, dir für $z$ gefunden wurden, voneinander unterscheiden, wollen wir anstelle von $\zeta$ und $\eta$ die angenommenen Werte schreiben, natürlich
\[
	\zeta = -c-f ~~ \mathrm{und} ~~ \eta = -c-g
\]
und die beiden ersten Reihen für den Buchstaben $z$ werden
\begin{align*}
z\ebinom{-g-1}{c-1} =& \ebinom{-g-1}{c-1} + \ebinom{-c-f}{1}\ebinom{-g-1}{c}x \\
&+ \ebinom{-c-f}{2}\ebinom{-g-1}{c+1}x^2 + \mathrm{etc.} \\
z\ebinom{-f-1}{c-1} =& \ebinom{-f-1}{c-1} + \ebinom{-c-g}{1}\ebinom{-f-1}{c}x \\
&+ \ebinom{-c-g}{2}\ebinom{-f-1}{c+1}x^2 + \mathrm{etc.}
\end{align*}
sein.

\paragraph{§18}
Damit wird diese Reihen auf eine gefälligere Form bringen, wollen wir
\[
	g+c-1 = h ~~ \mathrm{und} ~~ c-1=e
\]
setzen, so dass
\[
	c=e+1 ~~ \mathrm{und} ~~ g=h-e
\]
ist; daher wird nämlich unsere anfängliche Reihe
\begin{align*}
s\ebinom{h}{e} =& \ebinom{h}{e} + \ebinom{f}{1}\ebinom{h}{e+1}x \\
&+ \ebinom{f}{2}\ebinom{h}{e+2}x^2 + \ebinom{f}{3}\ebinom{h}{e+3}x^3 + \mathrm{etc.}
\end{align*}
sein. Die beiden folgenden Reihen, die aus dem Buchstaben $z$ gebildet wurden, werden aber zuerst
\begin{align*}
z\ebinom{e-h-1}{e} =& \ebinom{e-h-1}{e} + \ebinom{-e-f-1}{1}\ebinom{e-h-1}{e+1}x \\
&+ \ebinom{-e-f-1}{2}\ebinom{e-h-1}{e+2}x^2 + \mathrm{etc.}
\end{align*}
und darauf
\begin{align*}
z\ebinom{-f-1}{e} =& \ebinom{-f-1}{e} + \ebinom{-1-h}{1}\ebinom{-f-1}{e+1}x \\
&+ \ebinom{-1-h}{2}\ebinom{-f-1}{e+2}x^2 + \mathrm{etc.}
\end{align*}
sein; beide Größen $s$ und $z$ aber hängen voneinander so ab, dass
\[
	s = (1-x)^{f+b+1}z
\]
ist.

\paragraph{§19}
Wir wollen den immensen Nutzen der Transformation bei einem aus der Integralformel
\[
\int \frac{\partial \phi \cos{i\phi}}{(1+aa-2a\cos{\phi})^{n+1}}
\]
entstandenen höchst bemerkenswerten Fall zeigen, deren Integral von der Grenze $\phi = 0^{\circ}$ bis hin zur Grenze $\phi = 180^{\circ}$ ich freilich zuerst durch eine Vermutung allein geschlossen habe gleich
\[
	\frac{\pi a^i}{(1-aa)^{2n+1}}V
\]
zu sein, während
\begin{align*}
V = \ebinom{n-1}{0}\ebinom{n+i}{i} + \ebinom{n-i}{1}\ebinom{n+i}{i+1}aa + \ebinom{n-i}{2}\ebinom{n+i}{i+2}a^4 + \mathrm{etc.}
\end{align*}
wird; diese Reihe wird, wenn sie mit unserer anfänglichen zusammengebracht wird, dass $V = s\ebinom{h}{e}$ ist,
\[
	h = n+i ~~ \mathrm{und} ~~ e=i
\]
liefern, dann aber
\[
	f = n-i ~~ \mathrm{und} ~~ x = aa
\]
Die beiden anderen daher gebildeten Reihen werden zuerst
\begin{align*}
z\ebinom{-n-1}{i} =& \ebinom{-n-1}{i} + \ebinom{-n-1}{i}\ebinom{-n-1}{i+1}a^2 \\
&+ \ebinom{-n-1}{2}\ebinom{-n-1}{i+2}a^4 + \mathrm{etc.}
\end{align*}
andererseits
\begin{align*}
z\ebinom{i-n-1}{i} =& \ebinom{i-n-1}{i} + \ebinom{-n-i-1}{1}\ebinom{i-n-1}{i+1}a^2 \\
&+ \ebinom{-n-i-1}{2}\ebinom{i-n-1}{i+2}a^4 + \mathrm{etc.}
\end{align*}
sein, welche Reihen aus der Reihe $V$ selbst entstehen, indem man anstelle von $n$ gleich $-n-1$ schreibt. Nun wird aber die Beziehung zwischen $s$ und $z$
\[
	s = (1-aa)^{2n+1}z
\]
sein; dann ist aber
\[
	V = s\ebinom{n+i}{i}
\]

\paragraph{§20}
Weil daher
\[
\int\frac{\partial \phi \cos{i\phi}}{(1+aa-2a\cos{\phi})^{n+1}} = \frac{\pi a^i}{(1-aa)^{2n+1}}V = \frac{\pi a^i}{(1-aa)^{2n+1}}\ebinom{n+i}{i}s
\]
ist, wollen wir in dieser Form anstelle von $n$ gleich $-n-1$ schreiben und es sei
\[
\int\frac{\partial \phi \cos{i\phi}}{(1+aa-2a\cos{\phi})^{-n}} ~ \begin{bmatrix} \text{von } \phi = 0^{\circ} \\ \text{bis } \phi = 180^{\circ} \end{bmatrix} = \frac{\pi a^i}{(1-aa)^{-2n-1}}U
\]
wird
\[
	U = \ebinom{-n-1-i}{0}\ebinom{-n-1+i}{i} + \ebinom{-n-1-i}{1}\ebinom{-n-1+i}{i+1}aa + \mathrm{etc.}
\]
sein und daher
\[
	U = z\ebinom{i-n-1}{i} = \ebinom{i-n-1}{i}(1-aa)^{-2n-1}s
\]

\paragraph{§21}
Wir wollen gleich
\[
	1 + aa - 2a\cos{\phi} = \Delta
\]
setzen und diese zwei Integralwerte betrachten, die wir gerade erreicht haben:
\begin{align*}
\mathrm{I.}  \int\frac{\partial\phi\cos{i\phi}}{\Delta^{n+1}} &= \frac{\pi a^i}{(1-aa)^{2n+1}}\ebinom{n+i}{i}s \\
\mathrm{II.}  \int \Delta^n\partial\phi\cos{i\phi} &= \frac{\pi a^i}{(1-aa)^{-2n-1}}\ebinom{i-n-1}{i}(1-aa)^{-2n-1}s = \pi a^i \ebinom{i-n-1}{i}s
\end{align*}
Als logische Konsequenz folgern wir zwischen diesen zwei Integralformeln, die von der Grenze $\phi = 0^{\circ}$ bis zur Grenze $\phi = 180^{\circ}$ erstreckt wurden, diese höchst bemerkenswerte Relation:
\[
	\int\frac{\partial\phi\cos{i\phi}}{\Delta^{n+1}} :  \int \Delta^n\partial\cos{i\phi} = \ebinom{n+i}{i} : \ebinom{i-n-1}{i}(1-aa)^{2n+1}
\]
oder es wird
\[
	\ebinom{n+i}{i}(1-aa)^{-n}\int\Delta^n\partial\phi\cos{i\phi} = \ebinom{-n-1+i}{i}(1-aa)^{n+1}\int\Delta^{-n-1}\partial\phi\cos{i\phi}
\]
sein.

\paragraph{§22}
Ich hatte dieses letzte Theorem schon vor einiger Zeit auch durch Induktion allein gefunden, und ich hatte von einem Beweis davon schon fast die Hoffnung verloren, welcher nun aus der erwähnten Transformation der Reihen sich quasi wie von selbst offenbart; daher erkennt man den Nutzen dieser Transformation, die mit Recht von sehr weitem Umfang anzusehen ist, umso klarer.

\paragraph{§23}
Nachdem ich aber neulich dasselbe Theorem vorgelegt hatte, scheint sich das dort gegebene Verhältnis zwischen den beiden Integralformeln irgendwie von dem hier gefundenen zu unterscheiden; dennoch entdeckt man, dass sie perfekt übereinstimmen, wenn nur das folgende Verhältnis zur Hilfe genommen wird, wodurch im Allgemeinen
\[
	\ebinom{n}{i} : \ebinom{-n-1}{i} = \ebinom{-n-1+i}{i} : \ebinom{n+i}{i}
\]
ist, wessen Verhältnis daher natürlich klar ist, weil im Allgemeinen immer
\[
	\ebinom{-a}{i} = \pm\ebinom{a+i-1}{i}
\]
ist und daher auch 
\[
	\ebinom{b}{i} = \pm\ebinom{b-1+i}{i}
\]
wo die oberen Zeichen gelten, wenn $i$ eine gerade Zahl war, die unteren aber wenn ungerade. Daher wird also
\[
	\ebinom{n+i}{i} = \pm\ebinom{-n-1}{i} ~~ \mathrm{und} ~~ \ebinom{-n-1+i}{i} = \pm\ebinom{n}{i}
\]
sein.

\paragraph{§24}
Daher kann also unser Theorem noch gefälliger ausgesprochen werden. Wenn wir der Kürze wegen
\[
	\frac{1+aa-2a\cos{\phi}}{1-aa} = \Theta
\]
setzen, so dass
\[
	\Delta = (1-aa)\Theta
\]
ist, dann wird dieses Verhältnis hervorgehen:
\begin{align*}
	&\int\Theta^n\partial\phi\cos{i\phi} : \int\frac{\partial\phi\cos{i\phi}}{\Theta^{n+1}} = \int\frac{\Delta^n\partial\cos{i\phi}}{(1-aa)^n} : \int\frac{\partial\phi\cos{i\phi}(1-aa)^{n+1}}{\Delta^{n+1}} \\
	&= \ebinom{-n-1+i}{i}:\ebinom{n+i}{i} = \ebinom{n}{i} : \ebinom{-n-1}{i}
\end{align*}
und so wird
\[
	\ebinom{n}{i}\int\frac{\partial\phi\cos{i\phi}}{\Theta^{n+1}} = \ebinom{-n-1}{i}\int\Theta^n\partial\phi\cos{i\phi}
\]
sein.

\section*{Theorem}
\paragraph{§25}
Wenn die Summe dieser Reihe bekannt war:
\[
\frac{h}{e} + \ebinom{f}{1}\ebinom{h}{e+1}x + \ebinom{f}{2}\ebinom{h}{e+2}x^2 + \ebinom{f}{3}\ebinom{h}{e+3}x^3 + \mathrm{etc.}
\]
die man gleich $s$ setze, dann können auch die Summen der beiden folgenden Reihen beschafft werden, von denen die erste diese ist:
\[
\ebinom{e-h-1}{e} + \ebinom{-e-f-1}{1}\ebinom{e-h-1}{e+1}x + \ebinom{-e-f-1}{2}\ebinom{e-h-1}{e+2}x^2 + \mathrm{etc.}
\]
deren Summe
\[
	\ebinom{e-h-1}{e}\frac{s}{\ebinom{h}{e}(1-x)^{f+h+1}}
\]
sein wird, wo man bemerke, dass
\[
	\ebinom{e-h-1}{e} = \pm \ebinom{h}{e}
\]
ist, wo das obere Zeichen gilt, wenn $i$ eine gerade Zahl war, das untere wenn eine ungerade; daher ist die Summe dieser Reihe
\[
	\frac{\pm s}{(1-x)^{f+h+1}}
\]
Die andere Reihe, deren Summe daher bestimmt werden kann, wird
\[
	\ebinom{-f-1}{e} + \ebinom{-h-1}{1}\ebinom{-f-1}{e+1}x + \ebinom{-h-1}{2}\ebinom{-f-1}{e+2}x^2 + \mathrm{etc.}
\]
sein, deren Summe
\[
	\ebinom{-f-1}{e}\frac{s}{\ebinom{h}{e}(1-x)^{f+h+1}}
\]
sein wird, die auch auf diese Weise ausgedrückt werden kann:
\[
	\pm \ebinom{f+e}{e}\frac{s}{\ebinom{h}{e}(1-x)^{f+h+1}}
\]

\paragraph{§26}
Wenn die Summe dieser drei Reihen gesetzt werden wie folgt
\[
	\mathfrak{A} = \ebinom{h}{e} + \ebinom{f}{1}\ebinom{h}{e+1}x + \ebinom{f}{2}\ebinom{h}{e+2}x^2 + \ebinom{f}{3}\ebinom{h}{e+3}x^3 + \mathrm{etc.}
\]
\[
	\mathfrak{B} = \ebinom{e-h-1}{e} + \ebinom{-e-f-1}{1}\ebinom{e-h-1}{e+1}x + \ebinom{-e-f-1}{2}\ebinom{e-h-1}{e+2}x^2 + \mathrm{etc.}
\]
\[
	\mathfrak{C} = \ebinom{-f-1}{e} + \ebinom{-h-1}{1}\ebinom{-f-1}{e+1}x + \ebinom{-h-1}{2}\ebinom{-f-1}{e+2}x^2 + \mathrm{etc.}
\]
verhalten sie sich untereinander so, dass
\[
	\ebinom{e-h-1}{e}\mathfrak{A} = \ebinom{h}{e}(1-x)^{f+h+1}\mathfrak{B}
\]
\[
	\ebinom{-f-1}{e}\mathfrak{A} = \ebinom{h}{e}(1-x)^{f+h+1}\mathfrak{C}
\]
\[
	\ebinom{-f-1}{e}\mathfrak{B} = \ebinom{e-h-1}{e}\mathfrak{C}
\]
ist.
\end{document}